\documentclass[11pt]{article}
\usepackage{xypic}

\font\tenmsb=msbm10
\font\sevenmsb=msbm7
\font\fivemsb=msbm5

\newfam\msbfam
\textfont\msbfam=\tenmsb
\scriptfont\msbfam=\sevenmsb
\scriptscriptfont\msbfam=\fivemsb
\def\Bbb#1{{\fam\msbfam #1}}

\newcommand\qed{{\hspace*{\fill}Q.E.D.\vskip12pt plus 1pt}}

\newcommand\sA{{\cal A}}

\newcommand\sE{{\cal E}}

\newcommand\sH{{\cal H}}
\newcommand\sI{{\cal I}}

\newcommand\sL{{\cal L}}

\newcommand\sO{{\cal O}}

\newcommand\sQ{{\cal Q}}

\newcommand\sZ{{\cal Z}}

\newcommand\zed{{\Bbb Z}}
\newcommand\rat{{\Bbb Q}}

\newcommand\comp{{\Bbb C}}
\newcommand\hirz[1]{{\Bbb F}_{#1}}

\newcommand\pn[1]{{\Bbb P}^{#1}}
\newcommand\pnpair[2]{({\Bbb P}^{#1},\sO_{{\Bbb P}^{#1}}({#2}))}
\newcommand\pnsheaf[2]{\sO_{{\Bbb P}^{#1}}({#2})}
\newcommand\proj[1]{{\mathop{{\Bbb P}({#1})}\nolimits}}

\def\rank{{\mathop{\rm rank\,}\nolimits}}

\newcommand\proof{{\noindent\bf Proof.\ }}
\newcommand\Case{{\noindent\bf Case\ }}

\newtheorem{theorem}{Theorem}[section]
\newtheorem{lemma}[theorem]{Lemma}
\newtheorem{corollary}[theorem]{Corollary}
\newtheorem{proposition}[theorem]{Proposition}

\newtheorem{re}[theorem]{Remark}

\newenvironment{remark}{\begin{re}\em}{\end{re}}

\begin{document}
\title {Chern Numbers of Ample  Vector Bundles on Toric Surfaces}
\author{Sandra Di Rocco\and Andrew J. Sommese}
\date{\ \ }

\maketitle

\section*{Introduction} Let $\sE$ be an ample rank $r$ bundle on a
smooth  toric projective surface, $S$, whose topological Euler
characteristic is $e(S)$. In this article, we prove a number of
surprisingly strong lower bounds for $c_1(\sE)^2$ and $c_2(\sE)$.

First, we  show Corollary (\ref{easyLowerBoundForC1Square}), which
says that,  given $S$ and $\sE$ as above, if $e(S)\ge 5$, then
$c_1(\sE)^2\ge r^2e(S)$. Though simple, this is much stronger
than the known lower bounds over not necessarily toric surfaces.
For example, see \cite[Lemma 2.2]{BSS94}, where it is shown that
there are many rank two ample vector bundles
with $(c_1(\sE)^2,c_2(\sE))=(2,1)$ on products of two smooth curves,
 at least one of which
has positive genus.

We then prove an estimate, Theorem \ref{degree}, which is quite
strong for large $e(S)$ and $r$. As $e(S)$ goes to $\infty$ with
$r$ fixed, the leading term of this lower bound is $\displaystyle
(4r+2)e(S)\ln_2(e(S)/12)$, while if $e$ is fixed and $r$ goes to
$\infty$,  the leading term of this lower bound is $\displaystyle
3(e(S)-4)r^2$.  For example, $c_1^2(\sE)\ge 3r^2e(S)$, for
  $r\le 3$ if  $e(S)\ge 13$, or for
  $r\le 6$ if  $e(S)\ge 19$, or for
  $r\le 141$ if  $e(S)\ge 100$.  Or again,
  $c_1^2(\sE)\ge 5r^2e(S)$, for
  $r\le 10$ if  $e(S)\ge 100$.
   We include a three line
Maple program in Remark \ref{mapleProgram} for plotting the expression
for the lower bound.

The strategy is to use the adjunction process to find lower bounds for
$c_1(\sE)^2$. Toric geometry has two major implications for the
adjunction process. First, given an ample rank $r$ vector bundle
$\sE$ on a smooth toric surface $S$, there is the inequality
$-\det \sE \cdot K_S\ge e(S)(\rank \sE) $.  Adjunction theory yields
the  lower bound for $c_1(\sE)^2$ given in Theorem
\ref{easyLowerBoundForC1Square}, which implies that
$\displaystyle c_1(\sE)^2
>r^2e(S)$ for $e(S)\ge 7$. The
second important fact is that $h^0(tK_S+\det \sE)>0$ for integers $t$
between $0$ and at least $\rank \sE+\ln_2(e(S)/6)$.   Adjunction theory
yields the strong lower bound given in Theorem (\ref{degree}) for
$\displaystyle c_1(\sE)^2$ when $e(S)\ge 7$.

Using Bogomolov's instability theorem, we get the strong lower
bound given in Theorem (\ref{applicationOfBogomolov}) for the
second Chern class,  $c_2(\sE)$, of a rank two ample vector
bundle.  Basically if $c_2(\sE)$ is less than one fourth the lower
bound already derived for $c_1(\sE)^2$, then we have an unstable
bundle, and Bogomolov's instability theorem combined with the
Hodge index theorem give strong enough conditions to get a contradiction.  The
short list of exceptions to the bound $c_2(\sE)>e(S)$ are
classified.  Even assuming $\sE$  very ample on a nontoric
surface, the best general result \cite{BSS96} shows only that
$c_2(\sE)\ge 1$ with equality for $\pn 2$.

Inequalities derived from adjunction theory usually have the form,
``some inequality is true if certain projective invariants are
large enough.'' Typically examples exist outside the range where
the  adjunction theoretic method works.  For rank two ample vector
bundles $\sE$ we use a variety of special methods, including
adjunction theory and Bogomolov's instability theorem, to
enumerate the exceptions to  either the inequality $c_1(\sE)^2\ge
4e(S)$ or the inequality  $c_2(\sE)\ge e(S)$ holding. The
exceptions are collected in Table  1.

We would like to thank the Department of Mathematics of the K.T.H.
(Royal Institute of Technology) of Stockholm, Sweden, for making
our collaboration possible.  The second author would like to thank
the Department of Mathematics of Colorado State University for
their support and fine working environment during the period when
the final research for this work was carried out.


\section{Background material}\label{backgroundSection}In this paper we
 work over
$\comp$. By a variety we mean a complex analytic space, which might be
neither reduced or irreducible.

A rank $2$ vector bundle $\sE$ on a nonsingular surface $S$ is called
{\em Bogomolov unstable} \cite{R78}, or {\em unstable}  for short, if
$c_1(\sE)^2>4c_2(\sE)$. When $\sE$ is unstable there exists a line bundle
$\sA$ and a zero subscheme $(\sZ,\sO_\sZ)$ fitting in the exact sequence
 \begin{equation}\label{BS}
 0\to \sA\to \sE\to (\det\sE-\sA)\otimes\sI_\sZ\to 0;
 \end{equation}
with the property  that for all ample line bundles $\sL$ on $S$,
 $(2\sA-\det \sE)\cdot \sL >0$.  The standard consequences of this
result that we will often use in this article are:
 \begin{enumerate}
 \item $(2\sA-\det \sE)\cdot (2\sA-\det \sE)>0$, and $2\sA-\det \sE$ is
 $\rat$-effective;
 and
 \item for all nef and big line bundles $\sL$ on $S$,
  $(2\sA-\det \sE)\cdot \sL >0$.
 \end{enumerate}
We define $\sH:=\det \sE$. Note that
\begin{itemize}
\item $c_2(\sE)=\sA\cdot(\sH-\sA) +\deg(\sZ)$, where
$\deg(\sZ)=h^0(\sO_\sZ)$; and
\item the line bundle $\sH-\sA$ is a quotient of $\sE$
off a codimension two subset and therefore it is ample when $\sE$ is
ample.
\end{itemize}
Using the Hodge inequality
 $(\sH-\sA)^2\ (2\sA-\sH)^2\leq \left[(\sH-\sA)\cdot
(2\sA-\sH)\right]^2$, we obtain the following:
\begin{equation}\label{EQ1}
\sA\cdot(\sH-\sA)\geq (\sH-\sA)^2+\sqrt{(\sH-\sA)^2}
\end{equation}

A  toric surface $S$ is a
surface containing a two dimensional torus as Zariski open subset
and such that the action of the torus on itself extends to $S$. All
toric surfaces are normal. In this article we consider surfaces
polarized by an ample vector bundle, therefore $S$ will always denote
a normal projective toric surface.
For basic definitions on toric varieties we refer to \cite{O88}.

We recall that if $e:=e(S)$ is the Euler characteristics of $S$ then \
$\rank({\rm Pic}(S))=e-2\ $ and
    $\ K_{S}^{2}=12-e$.

We need the following useful lemmas, which are probably well known.

\begin{lemma}\label{vb}
    Let $\sE$ be a vector bundle over a normal $n$-dimensional toric variety.
   Assume $\proj{\sE}$ is toric, then $\sE=\oplus L_{i}$ with $L_{i}$
    equivariant line bundles.

\end{lemma}
    \proof Consider the bundle map $\proj\sE\to X$ with fiber
  $F=\pn{r-1}$ where $r:=rank(\sE)$. Every fiber has $r$-fixed
    points which define an unramified $r$ to one cover of $X$,
  $p:Y\to X$. $X$ being a normal toric variety, and thus simply connected,
     implies $Y=\cup
    X_{i}$ and $\sE=\oplus L_{i}$. \qed

It is classical  \cite{L82} that a surjective
  morphism $p : X\to Y$, with connected fibers between
  normal projective varieties, induces a
  homomorphism, from the connected
  component of the identity of the automorphism group of  $X$ to
   the connected component  of the identity of the automorphism group of
   $Y$, with respect to which $p$ is equivariant.  Using this basic fact
    we have the following lemma.
   \begin{lemma}\label{map}
   Let $p:X\to Y$ a surjective morphism with connected fibers
   from a   normal toric variety $X$ onto a normal variety $Y$.
   Then $Y$ and the general fiber of $p$ are toric.
   \end{lemma}

\begin{corollary}\label{singFibers}
Let $L$ be an ample line bundle on a smooth  projective toric
surface $S$.  If $f: S\to \pn 1$ is a morphism with connected
fibers, then the general fiber $F$ is isomorphic to $\pn 1$, there
are at most two singular fibers, and $e(S)\le 2+2L\cdot F$.
\end{corollary}
\proof Since the general fiber is toric it is isomorphic to $\pn 1$. From
equivariance we see that any singular fiber must lie over the two fixed
points of $\pn 1$.   Since there are at most $L\cdot F$ irreducible
components in a fiber, and there are at most two singular fibers the
inequality follows by considering the cases of no, one, or two singular
fibers. \qed

\begin{corollary}\label{simpleBlowup}
Let $f: S\to S'$ express a smooth toric surface $S$ as the blowup of a
smooth projective surface $S'$ at a finite set $B$.  Then $e(S)\le
2e(S')$.
\end{corollary}
\proof Let $b:= e(B)$, i.e., $b$ equals the cardinality of the finite set
$B$. Then we have $e(S)=e(S')+b$.  Since $S'$ is toric and $B$ are fixed
points of the toric action, we conclude that $e(B)$ is bounded by the
cardinality of the set of toric fixed points on $S'$, which is equal the
Euler characteristic of $S'$.  Thus we have $e(S)=e(S')+b\le 2e(S')$.
\qed

 Let $S$ be an irreducible toric surface. Then under the
prescribed torus action there are $e:=e(S)$ one dimensional
orbits. Denote their closures  by $D_i$ where $1\le i\le e$.  We
have the fundamental fact that
\begin{equation}\label{canonicalBundleFormula}
-K_S=\sum_{i=1}^{e(S)}D_i.
\end{equation}

We begin with a very simple observation which is in fact an
important tool in all our main results:

\begin{lemma}\label{KL}
    Let $\sE$ be an ample rank $r$ vector bundle on a projective
    normal toric surface
    $S$, and let  $\sH$ denote $\det\sE$. Then $-K_{S}\cdot\sH\geq re(S)$.
 \end{lemma}

\proof Let $\sH:=\det\sE=\sum_{1}^{e}a_{i} D_{i}$. By ampleness
$\sH\cdot D_{i}\geq r$ for all $i=1,\ldots,e$. Since
$K_{S}=\sum_{1}^{e}(-D_{i})$ we have
 $\displaystyle
 -K_{S}\cdot \sH=\sum_{1}^{e}\sH\cdot D_{i}\geq er.
 $
\qed

\begin{remark} In order to obtain the results in this paper we
use the bound \ref{KL} for $-KL$. The following example shows that
in general we cannot hope for a  better bound then the above.

Consider the toric surface given by the fan below, spanned by $12$
edges $\{\rho_{i}\}$ and with $12$ $2$-cones, i.e., $12$ fixed
points. The number before each edge indicates the self intersection
of the associated invariant divisor $D_{i}$.

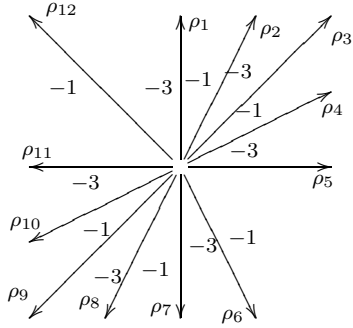
\[\begin{array}{ccccc}
 \xymatrix{
  & & & &  \\
 & & & & \\
 & & \uuto^(0.5){-3}_(0.9){\rho_{1}}\uurrto^(0.5){-3}_(0.9){\rho_{3}}
 \ar^(0.5)
 {-1}_(0.9){\rho_{2}}  @{->}[ruu] \ar^(0.5){-1}_(0.9){\rho_{4}} @{->}[rru]
 \rrto^(0.4){-3}_(0.9){\rho_{5}}\ar^(0.5){-1}_(0.9){\rho_{6}}@{->}[rdd]\ddto^(0.5)
 {-3}_(0.9){\rho_{7}}
 \ar^(0.6){-1}_(0.9){\rho_{8}} @{->}[ldd]\ddllto^(0.6)
 {-3}_(0.9){\rho_{9}}\ar^(0.6){-1}_(0.9){\rho_{10}} @{->}[lld]\llto^(0.6){-3}_(0.9)
 {\rho_{11}}
 \uullto^(0.6){-1}_(0.9){\rho_{12}} & & \\
  \\
  & & &
}
\end{array}\]

This surface is the equivariant blow up of $\pn{2}$ in $9$ points
and thus the Euler characteristics $e(S)=12$. Consider the line
bundle:
$$L=3D_{1}+5D_{2}+3D_{3}+5D_{4}+3D_{5}+5D_{6}+3D_{7}+5D_{8}+3D_{9}+3D_{10}+3D_{11}+
5D_{12}$$

It is ample since $L\cdot D_{i}=5-9+5=1$ for $i=1,3,5,7,9,11$ and $L\cdot
D_{i}=3-5+3=1$ for $i=2,4,6,8,10,12$. This also gives $\displaystyle
-LK_{S}=\sum_{1}^{12}L\cdot D_{i}=12=e$. Clearly this example can be
generalized to higher values of $e$.
\end{remark}

We end with a simple corollary of Lemma \ref{KL}.

\begin{corollary}\label{rSquareEvectorBundleLowerBound}
 Let $\sE$ be an ample rank $r$ vector bundle on a smooth
 projective toric surface $S$, and let $c_1^2:=c_1(\sE)^2$.
  If  $c_1^2\le re(S)$, then $r\le 3$ and either
  $g(\det \sE)=0$, and  $(S,\sE)$ is
 \begin{enumerate}
 \item   $\pnpair 2 1$ with $(c_1^2,e)=(1,3)$; or
 \item   $(\hirz 0,aE+bf)$ with $1\le ab\le 2$ and $(c_1^2,e)=(2ab,4)$; or
 \item $(\hirz 1,E+2f)$ with $(c_1^2,e)=(3,4)$; or
 \item $(\hirz 2,E+3f)$ with $(c_1^2,e)=(4,4)$; or
 \item $(\pn 2,\pnsheaf 21\oplus\pnsheaf 21)$ with $(c_1^2,e)=(4,3)$; or
 \end{enumerate}
    or $g(L)=1$, and $(S,\sE)$ is
 \begin{enumerate}
 \item   $(S,-K_S)$ with $(c_1^2,e)=(6,6)$; or
 \item   $(\hirz 0,(E+f)\oplus (E+f))$ with  $(c_1^2,e)=(8,4)$; or
 \item $(\pn 2,\pnsheaf 21\oplus\pnsheaf 21\oplus\pnsheaf 21)$ with
       $(c_1^2,e)=(9,3)$.
 \end{enumerate}
\end{corollary}
\proof

Let $\sH:=\det\sE$. If $\sH^2\le re$, then from $K_S\cdot \sH\le -re$ we
conclude that $2g(\sH)-2=\sH^2+K_S\cdot \sH\le 0$, and thus that
$g(\sH)\le 1$.

 If $g(L)=0$ we know
from classification theory, e.g., \cite{BS95,F90}, that $S$ is $\pn 2$ or
$\hirz r$.  A simple calculation shows the listed examples are the only
ones possible.

If $g(L)=1$, then from classification theory, e.g., \cite{BS95,F90}, we
know that $(S,\sH)$ is either a scroll over an elliptic curve or a Del
Pezzo surface with $\sH=-K_S$.  Since $S$ is toric and therefore
rational, $S$ is Del Pezzo. \qed

\section{Vector bundles over $\pn 2$ and $\hirz\epsilon$}\label{examples}
In this section we describe all pairs $(S,\sE)$ where $\sE$ is an ample
rank two bundle on a $\pn 2$ or a Hirzebruch surface, with the property
that either $c_1(\sE)^2\le 4e(S)$ or $c_2(\sE)\le e(S)$. Later in the
paper it will be shown  that these are all of the examples of rank
$2$ ample vector bundles $\sE$ on smooth toric surfaces, $S$, with either
$c_1(\sE)^2\le 4e(s)$ or $c_2(\sE)\leq e(S)$. The following table
includes the various cases. We give the Chern classes and indicate
whether the bundle is Bogomolov unstable ($U$), stable ($S$) or it is a
boundary case, i.e., $c_{1}^{2}=4c_{2}$, ($B$).

\begin{table}[htb]\label{theTable}\caption{All pairs $(S,\sE)$, with $\sE$
an ample rank two vector bundle on a smooth toric projective surface $S$,
and with either $c_1(\sE)^2\le 4e(S)$ or $c_2(\sE)\le e(S)$. The only
class where we do not know existence and uniqueness is listed on the last
line of the table.}
\begin{center}
\begin{tabular}{|c|c|c|c|c|c|} \hline
 $S $ &e(S)& $\sE$& $c_1(\sE)^2$& $c_2(\sE)$&$U/S/B$\\ \hline\hline
 $\pn{2}$&3&$\pnsheaf{2}1\oplus\pnsheaf{2}1$&$4$&$1$&$B$\\ \hline
 $\pn{2}$&3&$\pnsheaf{2}1\oplus\pnsheaf{2}2$&$9$&$2$&$U$\\ \hline
 $\pn{2}$&3&$T_{\pn{2}}$&$9$&$3$&$S$\\ \hline
 $\pn{2}$&3&$\pnsheaf{2}1\oplus\pnsheaf{2}3$&$16$&$3$&$U$\\ \hline
 $\pn{1}\times\pn{1}$&4&$p^*(\pnsheaf{1}1\oplus\pnsheaf{1}1)\otimes\xi$
  &$8$&$2$&$B$\\ \hline
 $\pn{1}\times\pn{1}$&4&$p^*(\pnsheaf{1}1\oplus\pnsheaf{1}2)\otimes\xi$&
  $12$&$3$&$B$\\ \hline
 $\pn{1}\times\pn{1}$&4&$p^*(\pnsheaf{1}1\oplus\pnsheaf{1}3)\otimes\xi$&
  $16$&$4$&$B$\\ \hline
 $\pn{1}\times\pn{1}$&4&$p^*(\pnsheaf{1}2\oplus\pnsheaf{1}2)\otimes\xi$&
   $16$&$4$&$B$\\ \hline
 $\pn{1}\times\pn{1}$&4&$\sO_{\pn{1}\times\pn{1}}(1,1)\oplus
   \sO_{\pn{1}\times\pn{1}}(2,2)$&$18$&$4$&$U$\\ \hline
 $\hirz{1}$&4& $p^*(\pnsheaf{1}1\oplus\pnsheaf{1}1)\otimes\xi$&
   $12$&$3$&$B$\\ \hline
 $\hirz{1}$&4&$p^*(\pnsheaf{1}1\oplus\pnsheaf{1}2)\otimes\xi$&
   $16$&$4$&$B$\\ \hline
 $\hirz{2}$&4&$p^*(\pnsheaf{1}1\oplus\pnsheaf{1}1)\otimes\xi$&
  $16$&$4$&$B$\\  \hline
 Del Pezzo&6& $(-K_S)\oplus(-K_S)$&24&6&$B$\\ \hline
  Del Pezzo&6&if any example exists, $\det\sE=-2K_S$ &24&$\ge 7$&$S$\\ \hline
\end{tabular}
\end{center}
\end{table}
 Fix the notation $c_2:=c_2(\sE)$, $\sH:=c_1=\det \sE$, and $e:=e(S)$. The
strategy that we follow is to first classify the pairs with
$c_1(\sE)^2\le 4e(S)$. Then any pair $(S,\sE)$ with $c_2\le e$ has
already been enumerated, or we have $c_2\le e < 4c_1^2$.  In the latter
case the bundle is unstable and we use the extra relations arising from
Bogomolov's instability theorem to classify the pair.

\subsection{$\pn{2}$}
Let $\sE$ be a rank two ample vector bundle over $\pn{2}$. Since
$\sH$ is the determinant bundle of a rank two bundle,
$\deg(\sH|_\ell)\geq 2$ for every line $\ell\in |\pnsheaf{2}1|$.
It follows that $\sH=\pnsheaf{2}a$ with $a\geq 2$. If $\sH^2\le
4e=12$, then $a=2,3$.  In case $a=2$, the restriction of $\sE$ to
each line $\ell$  is $\pnsheaf 11\oplus \pnsheaf 11$, and thus by
the classical results on uniform bundles \cite{OSS80}, $\sE=$\
\framebox{$\pnsheaf{2}1\oplus\pnsheaf{2}1$}. In case $a=3$, the
restriction of $\sE$ to each line $\ell$  is $\pnsheaf 11\oplus
\pnsheaf 12$, and thus by the classical results on uniform bundles
\cite{OSS80}, $\sE$=\ \framebox{$\pnsheaf{2}1\oplus\pnsheaf{2}2$}
or $\sE=$\ \framebox{$T_{\pn{2}}$}, the tangent bundle of $\pn 2$.

Now assume that $c_2(\sE)\leq 3$, but $c_1^2> 4e=12$. Thus it follows
that $\sH=\pnsheaf{2}a$ with $a\geq 4$. Since  $\sE$ is unstable, we have
a sequence as in (\ref{BS}) where $\sH-\sA= \pnsheaf{2}x$  and
$\sA=\pnsheaf{2}{ x+b}$ for $x,b>0$. The inequalities $3\geq
c_2(\sE)=x(x+b)+\deg(\sZ)$ and $a=2x+b\geq 4$ yield the only numerical
possibility: $(x,b+x)=(1,3)$ and $\deg(\sZ)=0$. Since $H^1(\pn
2,2\sA-\sH)=0$, we conclude the exact sequence splits, and it follows
that $\sE$=\ \framebox{$\pnsheaf{2}1\oplus\pnsheaf{2} 3$}.

\subsection{The Hirzebruch surfaces $\hirz{\epsilon}$}
Let $\hirz{\epsilon}=\proj{\sO_{\pn{1}}\oplus \pnsheaf{1}\epsilon}$ be the
Hirzebruch surface of degree $r$. Denote by $p:\proj{\sO_{\pn{1}}\otimes
\pnsheaf 1\epsilon}\to \pn{1}$ the projection map, and let  $F$ denote a
fiber of $p$. Let $\xi_\sE$ denote the tautological line bundle on
$\hirz{\epsilon}$, such that $p_*\xi_\sE\cong \sE$.
 Recall that ${\rm
Pic}(\hirz{r})=\zed F\oplus \zed E$, where $E$ is the section
corresponding to the surjection $\sO_{\pn{1}}\oplus \pnsheaf1\epsilon\to
\sO_{\pn{1}}$.  Note that $E^2=-\epsilon$.

The following is useful.
 \begin{lemma}\label{AmpleOnHirzRankR}
 Let $\sE$ be a rank $r$ ample vector bundle on $\hirz{\epsilon}$.
 Then $\det\sE\cdot F\ge r$ with equality if and only if
 $\sE\cong p^*V\otimes\xi_\sE$ where $V\cong\sE_E$.  In particular,
 in this case
 $$c_1(\sE)^2=r^2\epsilon+2r\det\sE\cdot E\ge r^2(2+\epsilon),$$
 and
 $$c_2(\sE)={r \choose 2}\epsilon+(r-1)\det\sE\cdot E\ge{r \choose 2}
   (2+\epsilon).$$
 \end{lemma}
\proof Since $\sE$ is a rank r ample vector bundle, and $F$ is a smooth
rational curve,
 we conclude that $\det\sE\cdot F\ge r$ with equality if and only if
 $\sE_F\cong \pnsheaf 1 1\oplus\cdots\oplus \pnsheaf 11.$
 In this case we have that $\sE\otimes \xi_\sE^*$ is trivial on every
 fiber and thus $\sE\otimes \xi_\sE^*\cong p^*V$ for some rank $r$ vector
 bundle on $\pn 1$.  Finally, note that $V\cong (p^*V)_E\cong \sE_E$.
 The rest of the lemma is a straightforward calculation.

 We record one simple corollary of the above Lemma.
 \begin{corollary}\label{rGreaterThan1}
 Let $\sE$ be a rank $r$ ample vector bundle on $\hirz{\epsilon}$.
 If $\epsilon \ge 2$ and $c_1(\sE)^2\le 4r^2$, then $\epsilon=2$ and
 $\sE\cong p^*(\pnsheaf 11\oplus\cdots\oplus\pnsheaf 11)\otimes\xi_\sE$.
 In this case $c_1(\sE)^2=4r^2$ and $c_2(\sE)=2r(r-1)$.
 \end{corollary}
 \proof  Let $\sH:=\det\sE =aE+bF$.
 Using Lemma \ref{AmpleOnHirzRankR}, we only need to show that
 $a=\sH\cdot F=r$.
 Assume therefore that $a\ge r+ 1$.  Then we have
 $\sH^2\ge a(2b-a\epsilon)\ge (r+1)(2r+(r+1)\epsilon)> 4r^2$. \qed

 Now assume that $c_2\le e=4$ or $c_1^2\le 4e=16$ and
$\sE|_{F}=\pnsheaf{1}a\oplus\pnsheaf{1}b$ with $a,b>0$.

\Case I: First consider the case when $(a,b)=(1,1)$. We are in the
situation of Lemma \ref{AmpleOnHirzRankR}. Letting
$V=\pnsheaf{1}\alpha\oplus\pnsheaf{1}\beta$, then
 $$4\geq c_{2}(\sE)=c_{2}(p^{*}(V)\otimes\xi)=\xi^{2}+\alpha+\beta
 =\epsilon+\alpha+\beta;
 $$
or
  $$
  12\geq c_1^2=c_{1}(p^{*}(V)\otimes\xi)^2=4\xi^{2}+4\alpha+4\beta
 =4(\epsilon+\alpha+\beta);
 $$
 The
only possible numerical possibilities are
 $\sE=$\ \framebox{$p^*(\pnsheaf{1}\alpha\oplus\pnsheaf{1}\beta)\otimes\xi$} with
$(\epsilon,\alpha,\beta)=(0,1,1),(0,1,2),(0,1,3),(0,2,2),
(1,1,1),(1,1,2),(2,1,1)$.

 \Case II: Assume now that $(a,b)\neq (1,1)$.  First,
  let us consider  the case $\epsilon=0$.
$\sH_F=\det(\sE)|_F=\pnsheaf{1}{a+b}$ implies $c_1^2\geq 18> 4e(S)$. Thus
if $c_2\le e=4$, $c_1^2\geq 4c_2(\sE)$ which means $\sE$ is unstable.
Consider the exact sequence (\ref{BS}).  We have that
$\sH-\sA=\sO_{\pn{1}\times\pn{1}}(x,y)$ for some $x>0$, $y> 0$, and
$\sA=\sO_{\pn{1}\times\pn{1}}(x+t,y+l)$ for some $t>0$, $l>0$. The
inequality $4\geq c_2(\sE)=x(y+l)+y(x+t)+\deg(\sZ)$ yields $\deg(\sZ)=0$
and $(x,y,x+t,y+l)=(1,1,2,2)$.  Since $\deg(\sZ)=0$ and
$H^1(\sO_{\pn{1}\times\pn{1}}(t,l))=0$,  we conclude that $\sE=$\
\framebox{$\sO_{\pn{1}\times\pn{1}}(1,1)\oplus
\sO_{\pn{1}\times\pn{1}}(2,2)$}\ .

Now assume that $\epsilon\geq 1$, and let $\sH=yF+xE$ with $x=a+b\geq 3$
and $\sH\cdot E=-x\epsilon+y\geq 2$ (since $\sH$ is the determinant of a
rank $2$ ample vector bundle). It follows that $\sH^2=a(2b-a\epsilon)\geq
a(4+a\epsilon)\geq 3(4+a\epsilon)\geq 21>4e(S)$. Thus if $c_2\le e=4$,
$c_1^2>4c_2(\sE)$ and thus $\sE$ is unstable. Let $\sA:=\alpha E_0+\beta
F$ be the line bundle in the sequence (\ref{BS}).   We have the following
straightforward inequalities:
 \begin{enumerate}
 \item $x\ge 3$, $\epsilon\ge 1$, $y\ge x\epsilon+2\ge 5$;
 \item $x-\alpha>0$, $y-\beta>0$, $y\ge \beta+(x-\alpha)\epsilon+1$;
 \item $2\alpha > x$, $2\beta>y$;
 \item $0<(2\sA-\sH)^2=(2\alpha-x)(4\beta-2y-(2\alpha-x)\epsilon)>0$,
 and in particular $4\beta+x\epsilon >2y+2\alpha \epsilon$; and
 \item $\sA\cdot (\sH-\sA)\le c_2(\sE)\le 4$, which gives
  $-\alpha(x-\alpha)\epsilon+\beta(x-\alpha)+\alpha(y-\beta)\le 4$.
 \end{enumerate}
Note that inequality (5) of the list can be written as
 $$
 \alpha(\alpha \epsilon-x-\beta+y) +\beta (x-\alpha)\le 4.
 $$
 Using inequality (2) from the list, $y-\beta\ge (x-\alpha)\epsilon +1$, we get
 \begin{eqnarray*}
 4&\ge& \alpha(\alpha \epsilon-x+(x-\alpha)\epsilon +1) +\beta(x-\alpha)\\
  &\ge& \alpha x (\epsilon-1) +\alpha +\beta(x-\alpha).
 \end{eqnarray*}
Now using equations (3) and (2) from the list we get the absurdity
 $$
  4 \ge  \alpha x(\epsilon-1) +\frac{x+1}{2} +\frac{y+1}{2}
    \ge      0      + \frac{4}{2} +\frac{5+1}{2}\ge 5.
 $$
 \qed

\section{Lower bounds for the Chern numbers of $\sE$}
In this section we obtain a number of  lower bounds for $c_1(\sE)^2$ for
a rank $r$ ample vector bundle on a smooth toric surface. Our main tool
is adjunction theory: good references for the standard adjunction results
that we use are \cite[Ch.\ 10, 11]{BS95} and \cite{F90}. The following is
a restatement, taking into account the geometry of toric surfaces, of the
main result for the adjunction theory for surfaces. Recall that on a
toric surface, a line bundle is ample if and only if it is very ample.

\begin{theorem}\label{adjunctionThm}
Let $L$ be an ample line bundle on a smooth  projective toric surface
$S$.
\begin{enumerate}
\item If $e=e(S)\geq 5$, then $K_S+L$ is spanned by global sections.
\item If $e=e(S)\geq 7$, then  $S$ is the equivariant blowup $\pi :
  S\to S_1$ of a smooth toric projective surface $S_1$ at a finite set $B$, such that
$L=\pi^{*}L'-\pi^{-1}(B)$ where $K_S+L\cong \pi^*(K_{S_1}+L')$, and both
$L'$ and $L_1:=K_{S_1}+L'$ are very ample.
\end{enumerate}
\end{theorem}
\proof Using \cite[9.2.2]{BS95},  note that the exceptions to $K_S+L$
being spanned by global sections are all ruled out by $e(S)\ge 5$.  The
associated map $p_{K_{S}+L}$ has a Remmert-Stein factorization $p=s\circ
\pi$ where $\pi:S\to S_{1}$ has connected fibers.  By Lemma (\ref{KL}),
we see that $e\ge 7$ rules out $\dim S_1=0$.  If $\dim S_1=1$, then we
have that $L\cdot F=2$ for a general fiber of $r$, but this  and $e\ge 7$
contradicts Corollary \ref{singFibers}.

 Since $\dim S_1=2$, it follows
from adjunction theory that $\pi : S\to S_1$ is the blow up of a smooth toric projective
surface $S_1$ at a finite set $B$, such that $L=\pi^{*}L'-\pi^{-1}(B)$
where $K_S+L\cong \pi^*(K_{S_1}+L')$, and both $L'$  and
$L_1:=K_{S_1}+L'$ are ample. The very ampleness of the last two bundles
follows from the fact that ample line bundles are  very ample on toric
varieties. \qed

 \begin{corollary}\label{easyLowerBoundForC1Square}
     Let $\sE$ be an ample rank $r$ vector
     bundle on a nonsingular toric surface $S$.
     If $e(S)\ge 5$ then $$c_1(\sE)^{2}\ge r^2e(S)$$
     with equality only if $\det\sE=-rK_S$ and $e(S)=6$.
 \end{corollary}
\proof Let $\sH:= \det \sE$. Let $t$ be the smallest positive integer for
which $tK_S+\sH$  is not ample.  Since $e(S)\ge 5$, $E\cdot (tK_S+\sH)=0$ for a
smooth rational curve $E$ with self-intersection $-1$. Thus
we have
 $$-t+E\cdot\sH=E\cdot (tK_S+\sH)= 0.$$
 Since $\sE$ has rank $r$, we have that $r\le\sH\cdot E=t$.  Thus
 $rK_S+\sH$ is spanned.  Using  Lemma \ref{KL},we have
 $$\sH^2\ge -\sH\cdot rK_S\ge r^2e(S).$$
 Moreover, since $\sH$ is ample, we have equality only if $\sH\cong
 -rK_S$.  In this case we have $r^2K_S^2=\sH^2= r^2e(S)$, or $K_S^2=
 e(S)$.  Since $K_S^2+e(S)=12$ we conclude that $K_S^2= 6$.
\qed

\begin{lemma}\label{delPezzo}     Let $\sE$ be an ample rank two vector
     bundle on a nonsingular toric surface $S$.  If $\det\sE=-2K_S$,
     $e(S)=6$, and $c_2(\sE)\le 6$, then $\sE:=-K_S\oplus
     -K_S$.
\end{lemma}
 \proof A simple computation shows that the Chern character of
 $\sE\otimes K_S$ is $2+(K_S^2-c_2(\sE))=2$.  Thus $\chi(\sE\otimes
 K_S)=2$.  Since $H^2(\sE\otimes K_S)=H^0(\sE^*)=0$, we conclude that
 $\dim H^0(\sE\otimes K_S)\ge 2$. Choose linearly independent $s_{1}, s_{2}\in
H^0(\sE\otimes K_S)$.

If $s_{1}\wedge s_{2}\neq 0$ then, since $\det(\sE\otimes K_S)=\sO_S$, we
conclude that $\sE\otimes K_S=\sO_S\oplus\sO_S$ i.e., $\sE\cong
-K_S\oplus -K_S$.

Thus we can assume without loss of generality that $s_{1}\wedge s_{2}=
0$.  The saturation $\sA$ of the images of $\sO_S$ in $\sE$, under
the two maps $g\to g\cdot s_i$, are equal.  $\sA$ is invertible, and
tensoring with $-K_S$ we have an exact sequence
 $$
 0\to \sA-K_S\to \sE \to \sQ \otimes\sI_\sZ\to 0,
 $$
 with $\sZ$ a $0$-dimensional subscheme of $S$.
 Note that $\sQ$ is ample, and therefore since $S$ is toric, very ample.
Since $e(S)=6$, we know that $S$ is not $\pn 2$ or a quadric, and thus
 \begin{equation}\label{sQsQLowerBound}
 \sQ^2\ge 3.
 \end{equation}
 Thus the  Hodge index theorem gives
 $(\sQ\cdot (-K_S))^2\ge \sQ^2(-K_S)^2\ge 18$, which implies that
 \begin{equation}\label{sQ-K_SLowerBound}
 \sQ\cdot (-K_S)\ge 5.
 \end{equation}
 Since $h^0(\sA)\ge 2$,
 we have $\sQ\cdot\sA\ge 1$.  Using this, and equations
 (\ref{sQsQLowerBound}) and (\ref{sQ-K_SLowerBound}) we have
 $$
 6=c_2(\sE)=(\sA-K_S)\cdot\sQ+\deg \sZ\ge 1+5+\deg\sZ.
 $$
 Thus $\deg\sZ=0$ and $\sA\cdot\sQ=1$. The exact sequence
 $$0\to\sA-K_S\to \sE\to \sQ \to 0$$
gives $-2K_S=c_{1}(\sE)=\sA+\sQ-K_S$ and $K_S+\sA+\sQ=\sO$. Thus
$(K_S+\sQ)\cdot\sQ=-\sA\cdot\sQ=-1$.  This is absurd, since on any smooth
surface $S$, the parity of $(K_S+L)\cdot L$ is even for any line bundle
$L$. \qed

\begin{remark}
  We do not know if there are any examples of $\sE$ satisfying all the
  hypotheses of Lemma \ref{delPezzo}, except that $c_2(\sE)>6$.
\end{remark}

 \begin{remark}
 The only smooth toric surfaces $S$ with $e(S)\le 4$ are $\pn 2$ or
Hirzebruch surfaces. Corollary \ref{rSquareEvectorBundleLowerBound}
classifies the exceptions to $c_1^2(S)>r^2e(S)$ for $r=1$, and \S
\ref{examples} classifies the exceptions for $r=2$ and $e(S)\le 4$.  They
are contained in Table 1.  For $\pn 2$ it seems difficult to
classify the exceptions when $r\ge 3$.  For the Hirzebruch surfaces
$\hirz\epsilon$, Corollary \ref{rGreaterThan1} classifies the exceptions
if $\epsilon\ge 2$.
 \end{remark}

If $e(S_{1})\geq 7$, we can repeat the procedure in Theorem
\ref{adjunctionThm}, using $L_1$ on $S_1$ in the same way we used $L$ on
$S$, and get $(S_{2},L_{2})$. We say the procedure has terminated when we
reach the first integer $b$ with $e(S_{b})\le 6$. (See \cite{BL89} for a
further study of the adjunction process.)  We call the sequence
$(S,L),\dots, (S_{b}, L_{b})$ the iterated adjunction sequence and $b$
the adjunction length of $S$.

Notice that in the iterated adjunction sequence, at every step we
 contract down $(-1)$-lines in $S_{i}$ with respect to the
 polarization $K_{S_{i}}+L_{i}$. This implies by Corollary
 (\ref{simpleBlowup})
  that
 $e(S_{i+1})\geq \lfloor\frac{e(S_{i})}{2}\rfloor$. If we assume, to
 start with, that the surface $S$ has $e(S)\geq 2^{b-1}\cdot 6+1$ then the
 adjunction length is at least $b$.

We have the following strong bound.
 \begin{theorem}\label{degree}
     Let $S$ be a nonsingular toric surface with
     $2^{b}\cdot 12\ge e(S)\geq 2^{b}\cdot 6+1$ for
     some integer $b\ge 0$ and e:=e(S).
     Let $\sE$ be an ample rank $r$ vector
     bundle on $S$, then $$c_1(\sE)^{2}\geq
     e(3r^2+2r+4br+2b-2)-12(b+1)(b+2r)-12r(r-1)+\frac{e}{2^{b-1}}-2$$
     \end{theorem}

\proof Since $\sH:=\det(\sE)$ is the determinant of a rank $r$
ample vector bundle, there are no smooth rational curves $C$ on
the polarized surface $(S,\sH)$ with $\sH\cdot C\le r-1$.
Therefore by Theorem (\ref{adjunctionThm}), $L:=K+(r-1)\sH$ ample.
Using Lemma \ref{canonicalBundleFormula}, we have the bound
 \begin{equation}\label{estimate}
 -K\sH\geq re.
 \end{equation}

The assumption $e(S)\geq 2^{b}\cdot 6+1$ implies that we have the
adjunction sequence  $(S,L),\dots, (S_{b}, L_{b}), (S_{b+1},
L_{b+1})$ with $L_{b+1}$ very ample.  It follows that the
sectional genus $g(L_{b+1})=g(K_{S_{b}}+L_{b})\geq 0$, i.e.,
$(K_{S_{b}}+L_{b})\cdot (K_{S_{b}}+K_{S_{b}}+L_{b})\geq -2$.

Let $S\to S_{1}\to\ldots\to S_{b}$ the sequence of contractions and let
$\pi_{i}$ denote the $i$-th contraction map. For simplicity let us set
$K_{i} :=(\pi\circ \pi_{1}\ldots \circ \pi_{i})^{*}(K_{S_{i}})$, $K_0:=
K_S$, and $S:=S_0$.
 $$(K_{S_{b}}+L_{b})\cdot
 (K_{S_{b}}+K_{S_{b}}+L_{b})=(K_{b}+K_{b-1}+\ldots +K_{1}+K_0+L)\cdot
 (K_b+K_b+K_{b-1}+\ldots +K_{1}+K_0+L)$$
 We can further decompose:
 \begin{eqnarray*}
  K_{b}\cdot(K_b+K_{b-1}+\ldots +K_{1}+K_0+L)&=&K_b^2+K_b\cdot(K_{b-1}+
                                 K_{b-2}+\ldots +K_{1}+K_0+L)\\
                                      &=&K_b^2+K_{b-1}\cdot(K_{b-1}+
                                 K_{b-2}+\ldots +K_{1}+K_0+L)\\
                                      &=&K_b^2+K_{b-1}^2+K_{b-1}\cdot(K_{b-2}+
                                         \ldots +K_{1}+K_0+L)\\
                                      &\vdots&\vdots\\
  &=&K_b^2+K_{b-1}^{2}+K_{b-2}^{2}+\ldots +
           K_{1}^{2}+K_0^2+K_0\cdot L\\
  (K_b+K_{b-1}+K_{b-2}+\ldots+K_{1}+K_0+L)^{2}&=&K_b^2+
      2K_{b}\cdot(K_{b-1}+\ldots+K_{1}+K_0+L)\\
      &&+(K_{b-1}+\ldots +K_{1}+K_0+L)^{2}\\
     &=&K_b^2+2(K_{b-1}^2+K_{b-2}^2+\ldots+K_1^2+K_0^2+K_0\cdot L)\\
      &&K_{b-1}^2+2K_{b-1}\cdot(K_{b-2}+\ldots+K_1+K_0+L)\\
      &&+(K_{b-2}+\ldots +K_1+K_0+L)^{2}\\
      &=&K_b^2+3K_{b-1}^{2}+5K_{b-2}^{2}+7K_{b-3}^{2}+\ldots\\
  &&+(2b-1)K_1^{2}+(2b+1)K_0^2+(2b+2)K_0\cdot L+L^{2}
\end{eqnarray*}
 Then: \ \
 $\displaystyle (K_{S_{b}}+L_{b})\cdot
(K_{S_{b}}+K_{S_{b}}+L_{b}) =$
\begin{equation}\label{Deg}
2K_b^2+4K_{b-1}^{2}+6K_{b-2}^{2}+\ldots +2b K_{1}^2+(2b+2)
K_0^2+(2b+3)K_0\cdot L+L^{2}\geq -2
\end{equation}
Recall that $K_{i}^{2}=12-e(S_{i})$ and $
e(S_{i})\geq(\frac{e}{2^{i}})$. Then
 \begin{eqnarray*}
  L^{2}+(2b+3)K_0\cdot L&\geq& -2-2\left(12-\frac{e}{2^{b}}\right)
     -4\left(12-\frac{e}{2^{b-1}}\right)-\ldots-(2b+2)(12-e)+(2b+3)e\\
        &\geq& -2
        -12(b+1)(b+2)+\frac{2e}{2^{b}}\sum_{j=0}^{b}((j+1)2^{j}).
 \end{eqnarray*}
 Using $\displaystyle\sum_{j=0}^{b}((j+1)2^{j})=2^{b+1}b+1$
we have
 \begin{eqnarray*}
 L^{2}+(2b+3)K_0\cdot L&\geq& -2-12(b+1)(b+2)+4eb+\frac{e}{2^{b-1}}.
  \end{eqnarray*}
Recalling equation (\ref{estimate}) and the fact that
$L=(r-1)K_0+\sH$, we get
\begin{eqnarray*}
\sH^{2}&\geq&
-2-12(b+1)(b+2)+4eb+\frac{e}{2^{b-1}}+2(r-1)re+(r-1)^2(e-12)\\
&&+(2b+3)re+(2b+3)(r-1)(e-12)\\
 &=& e(3r^2+2r+4br+2b-2)-12(b+1)(b+2r)-12r(r-1)+\frac{e}{2^{b-1}}-2.
 \end{eqnarray*}
 \qed

\begin{remark}\label{mapleProgram}
To get a global feel for the bound, we have found it helpful to graph the
expression. We include a short Maple V Release 5.1 program to plot the expression
divided by part of the leading term.  Varying the range of the rank $r$ and
the Euler characteristic $e$, and of the exact variant of \verb+lowerBound+,
 the scaled expression for the lower bound  is
useful.

\begin{verbatim}
b := floor(ln[2]((e-1)/6));
lowerBound := (r,e) -> e*(3*r^2+2*r+4*b*r+2*b-2)-12*(b+1)*(b+2*r)
                                                -12*r*(r-1)+e/2^(b-1)-2;
plot3d(lowerBound(r,e)/(r*e*(3*r+4*b)),r=1..20,e=13..100,style=PATCH,axes=BOXED);
\end{verbatim}
\end{remark}

\begin{remark} It is easily checked that the expression in $e$ and $r$
occurring in the lower bound is an increasing function of $e$ and $r$
for $e\ge 7$, $r\ge 1$.  It is also easy to check using the above bound
that $c_1(\sE)^2\ge 2r^2e(S)$ if $e(S)\ge 12$, and
$c_1(\sE)^2\ge 3r^2e(S)$ if $e(S)\ge 6r+7$.

Theorem (\ref{degree}) gives a strong asymptotic lower
bound for $c_1^2$ as $e$ goes to $\infty$. For any fixed $c>0$,
there will only be a finite number of possible pairs $(c_1^2,e)$
of numerical invariants for ample vector bundles $\sE$ on smooth
toric surfaces $S$ with $L^2\le ce$. For example, $c_1^2\ge
2r^2e(S)$ as soon as $e(S)\ge 13$.  This suggests that enumerating
the pairs $(S,\sE)$ with $\sH^2\le cre(S)$, where $\sE$ is an
ample vector bundle on a smooth toric surface $S$, and small $c>1$
should be a tractable classification problem with a nice answer.
\end{remark}

\begin{theorem}\label{applicationOfBogomolov}
 Let $\sE$ be an ample rank two vector bundle on a nonsingular
  toric variety $S$ with $2^{b}\cdot 12\ge e(S)\geq 2^{b}\cdot 6+1$ for
     some integer $b\ge 0$ and $e:=e(S)$. Then
$$c_{2}(\sE)\geq
-3(b+2)(b+3)+\frac{5b+7}{2}e+\frac{e}{2^{b+1}}-\frac{1}{2}$$
\end{theorem}
\proof If the inequality is not satisfied then using Theorem
(\ref{degree}), $c_{1}(\sE)^{2}>4c_{2}(\sE)$, and thus the bundle
would be unstable. The exact sequence (\ref{BS}) and the
inequality (\ref{EQ1}) give $$c_{2}(\sE)\geq
(\sH-\sA)^{2}+\sqrt{(\sH-\sA)^{2}}$$ the  divisor $\sH-\sA$ is
ample and thus by
 Theorem (\ref{degree})
 $$-3(b+2)(b+3)+\frac{(5b+7)}{2}e+\frac{e}{2^{b+1}}-\frac{1}{2}>
c_{2}(\sE)\ge e(6b+3)-12(b+1)(b+2)+\frac{e}{2^{b-1}}-2+1$$
 which is
equivalent to \ $18b^2+42b+13-7eb+e-3e/2^b>0$, which is impossible.\qed

\begin{remark}
We expect that a generalization of Theorem \ref{applicationOfBogomolov}
to ample vector bundles of arbitrary rank $r$ is true.  Based on a strong
dose of optimism, we conjecture that if $\sE$ is an ample rank $r$ vector
bundle on a smooth toric projective surface $S$ with $2^{b}\cdot 12\ge
e(S)\geq 2^{b}\cdot 6+1$ for some integer $b\ge 0$, then
 $$
 c_2(\sE)\ge \frac{r-1}{2r}\left[e(S)(3r^2+2r+4br+2b-2)-12(b+1)(b+2r)-
                 12r(r-1)+\frac{e}{2^{b-1}}-2\right].
 $$
\end{remark}

We now turn to the special case of rank two bundles where the inequality
$\displaystyle c_{2}({\sE})> e(S)$ fails to be true.

\begin{lemma}\label{generalBogRestriction}
Let $\sE$ be an unstable  ample rank two vector bundle on a smooth toric
projective surface $S$. If $\sE$ is Bogomolov unstable and
 $c_2(\sE) \le e(S)+\sqrt{e(S)}$, then $S$ is either $\pn 2$ or $\hirz
 \epsilon$ with $\epsilon\le 2$.
 \end{lemma}
\proof Assume that $\sE$ is Bogomolov unstable. Consider the sequence
(\ref{BS}) and the inequality:
 $$e(S)+\sqrt{e(S)}\geq
 c_{2}(\sE)=\sA\cdot (\sH-\sA)+\deg(\sZ)\geq
(\sH-\sA)^2+\sqrt{(\sH-\sA)^{2}}$$
 We can then assume $(\sH-\sA)^{2}\le e$.
We now apply Theorem \ref{rSquareEvectorBundleLowerBound} to the ample
line bundle $\sH-\sA$.
 \qed

\begin{remark}
Let $\delta := \min\{L^2 | L \text{ an ample line bundle on } S\}$. The
 above argument implies that any ample vector bundle $\sE$ with
$c_2(\sE)< \delta+ \sqrt{\delta}$ is Bogomolov stable.
\end{remark}

\begin{corollary}\label{Bog}Let $\sE$ be an ample rank
 two vector bundle on
a smooth toric projective surface $S$. Assume that $c_{2}(\sE)\le e(S)$,
if $\sE$ is not Bogomolov Stable  then $(S,\sE)$ is contained the Table
1.
\end{corollary}
\proof Simply use Lemma \ref{generalBogRestriction} and the results for
$\pn 2$ and the Hirzebruch surfaces from \S \ref{examples}
 \qed

 \begin{proposition}\label{main}
 Let $\sE$ be an ample rank two vector bundle on a
 smooth projective toric surface $S$. If either $c_1(\sE)^2\le 4e(S)$ or
 $c_{2}(\sE) \le e(S)$, then $(S,\sE)$ is in the Table 1.
 \end{proposition}

\proof We can also assume that $S$ is neither $\pn 2$ or a Hirzebruch
surface by using the results of \S 2. Thus $e(S)\ge 4$.  Using Corollary
\ref{easyLowerBoundForC1Square} and Lemma \ref{delPezzo}, we can assume
without loss of generality that $c_1(\sE)^2>4e(S)$. If $c_2(\sE)\leq e$,
then we are in the situation of  Lemma \ref{Bog}. \qed

\vspace{1cm} \small
 \begin{tabular}{lcl}
 Sandra Di Rocco                  && Andrew J. Sommese \\
 Department of Mathematics        &&Department of Mathematics\\
 KTH                              && University of Notre Dame\\
 S-100 44 Stockholm, Sweden       &&Notre Dame, Indiana 46556, U.S.A,\\
 fax: Sweden + 46 8 7231788       &&fax: U.S.A. + 219--631-6579 \\
 sandra@math.kth.se               && sommese@nd.edu\\
 URL:{\tt www.math.kth.se/$\sim$sandra}{\ \ \ \ \ \ } &&URL: {\tt www.nd.edu/$\sim$sommese}\\
 \end{tabular}

\end{document}